\documentstyle[11pt,amscd,amssymb]{amsart}
\theoremstyle{plain}
 \newtheorem{thm}{Theorem}[section]
 \newtheorem{lem}[thm]{Lemma}
 \newtheorem{prop}[thm]{Proposition}
 \newtheorem{cor}[thm]{Corollary}
 
\theoremstyle{definition}
 \newtheorem{defn}{Definition}[section]
\theoremstyle{remark}
 \newtheorem{rem}{Remark}[section]
 \newtheorem{ex}{Example}

\newcommand{\Hom}{\operatorname{Hom}}
\newcommand{\codim}{\operatorname{codim}}
\newcommand{\im}{\operatorname{im}}

\newcommand{\rk}{\operatorname{rk}}

\newcommand{\NS}{\operatorname{NS}}

\newcommand{\Pic}{\operatorname{Pic}}
\newcommand{\ch}{\operatorname{ch}}
\newcommand{\td}{\operatorname{td}}
\newcommand{\Alb}{\operatorname{Alb}}
\newcommand{\Hilb}{\operatorname{Hilb}}

\newcommand{\Spec}{\operatorname{Spec}}

\newcommand{\Km}{\operatorname{Km}}

\font\b=cmr10 scaled \magstep5
\def\bigzerou{\smash{\lower1.7ex\hbox{\b 0}}}

\numberwithin{equation}{section}
\begin{document}

\title{
Albanese map of moduli of stable sheaves on abelian surfaces
}
\author{K\={o}ta Yoshioka}
 \address{Max Planck Institut f\"{ur} Mathematik,
Gottfried Claren Str. 26, D-53225 Bonn, Germany $\&$
\newline
Department of mathematics, Faculty of Science, Kobe University,
Kobe, 657, Japan}
\email{yoshioka@@mpim-bonn.mpg.de}
 \maketitle

\section{Introduction}

Let $X$ be a smooth projective surface defined over ${\Bbb C}$
and $H$ an ample line bundle on $X$.
If $K_X$ is trivial,
Mukai [M3] introduced a quite useful notion called 
Mukai lattice $(H^{ev}(X,{\Bbb Z}),\langle \;\;,\;\;\rangle)$,
where $H^{ev}(X,{\Bbb Z})=\oplus_i H^{2i}(X,{\Bbb Z})$.
For a coherent sheaf $E$ on $X$,
we can attach an element of $H^{ev}(X,{\Bbb Z})$
called Mukai vector $v(E):=\ch(E)\sqrt{\td_X}$,
where $\td_X$ is the Todd class of $X$.
We denote the moduli space of stable sheaves $E$
of $v(E)=v$ by $M_H(v)$.
If $H$ is general (i.e. it does not lie on walls [Y1])
and $v$ is primitive, then $M_H(v)$ is a smooth projective scheme.

If $X$ is a K3 surface and $v$ is primitive, then
$M_H(v)$ is extensively studied by many authors.
In particular, in many cases,  
$M_H(v)$ is an irreducible sympletic
manifold and the period of $M_H(v)$ is written down in
terms of Mukai lattice ([Mu3,5], [O], [Y4]).

In this paper, we shall treat the case where $X$ is an abelian surface. 
In [Y2], we studied $H^i(M_H(v),{\Bbb Z})$
$i=1,2$ under some assumptions on $v$.
We also constructed a morphism ${\frak a}:M_H(v) \to X \times \widehat{X}$
and proved that ${\frak a}$ is an albanese map
for $\langle v^2 \rangle \geq 2$,
where $\widehat{X}$ is the dual of $X$. 
In this paper, we shall consider the fiber of albanese map
under the same assumptions on $v$ in [Y2] (cf. Theorem \ref{thm:period}).
If $\langle v^2 \rangle =0$, then Mukai showed that
$M_H(v)$ is an abelian surface (see [Mu5, (5.13)]).
In this case, ${\frak a}$ is an immersion.
If $\langle v^2 \rangle=2$, then Mukai [Mu1] and the author
[Y2, Prop. 4.2] showed that ${\frak a}:M_H(v) \to X \times \widehat{X}$
is an isomorphism.
Hence we assume that $\langle v^2 \rangle \geq 4$. 
Let $K_H(v)$ be a fiber of ${\frak a}$.
Then $\dim K_H(v)=\langle v^2 \rangle-2$.
Hence if $\langle v^2 \rangle \geq 6$, then 
$\dim K_H(v) \geq 4$.
In this case, we get the following,
which is an analogous result to that for a K3 surface.

\begin{thm}\label{thm:period}
Let $X$ be an abelian surface.
Let $v=r+\xi+a \omega \in
H^{ev}(X,{\Bbb Z})$, $\xi \in H^2(X,{\Bbb Z})$
be a Mukai vector such that $r>0$,
$r+\xi$ is primitive
and $\langle v^2 \rangle \geq 6$,
where $\omega$ is the fundamental class of $X$.
Then for a general ample line bundle $H$,
$K_H(v)$ is an irreducible symplectic manifold and
\begin{equation}
\theta_v:v^{\perp} \to H^2(K_H(v),{\Bbb Z})
\end{equation}
is an isometry of Hodge structures.
\end{thm}
For the definition of $\theta_v$, see preliminaries. 
Our theorem shows that Mukai lattice for an abelian surface
is as important as that for a K3 surface. 
As an application of this theorem, we shall show that
for some $v$, $M_H(v)$ is not birationally equivalent
to $\widehat{Y} \times \Hilb_Y^n$ for any $Y$
( Example \ref{ex:ex1}).

In section 1, we collect some known facts which will be used in this paper.
Since the canonical bundle of $M_H(v)$ is trivial,
$M_H(v)$ has a Bogomolov decomposition.
We shall also construct a decomposition which will become
a Bogomolov decomposition for $M_H(v)$.

In section 2, we shall prove Theorem \ref{thm:period}.  
We shall first treat rank 1 case.
In this case, $K_H(v)$ is the generalized Kummer variety
$K_{n-1}$
constructed by Beauville [B],
where $\langle v^2 \rangle/2=n$.
Hence Theorem \ref{thm:period} follows from 
Beauville's description of $H^2(K_{n-1},{\Bbb Q})$ and
some computations.  
For higher rank cases, we shall use the same method as in [Y2].
More precisely,
we shall first treat the case where $X$ is a product of two elliptic curves.
In this case, we constructed a family of stable sheaves
$E$ of $v(E)=v$ in [Y2].
Then it induces a birational map
$\widehat{X} \times \Hilb_X^n \cdots \to M_H(v)$, and hence we get a birational
map from the generalized Kummer variety $K_{n-1}$
to $K_H(v)$, where $n=\langle v^2 \rangle/2$. 
  By the description of $\theta_v(x), x \in v^{\perp}$ in [Y2],
we get our theorem for this case.
For general cases, we shall use deformation arguments as in 
[G-H], [O] and [Y2].

In section 3, we shall treat the remaining case.
In this case, we shall prove that $K_H(v)$ is isomorphic to
a moduli space of stable sheaves on the Kummer surface associated to $X$. 

In appendix, we shall explain a more sophisticated method to prove
 Theorem \ref{thm:period} at least for $r \ne 2$.
In the K3 cases, we know that isometries of Mukai lattice
are quite useful to compute period of moduli spaces [Mu3], [Y4].
Hence it is also important to study isometries in our cases.
Fourier-Mukai transforms are good examples of isometries [Mu4].
Thanks to recent results of Bridgeland [Br] on Fourier-Mukai
transforms,
 we can replace our computations in [Y2] 
to a simple calculation (Proposition \ref{prop:comm})
for $\geq 3$.
Thus we get another proof of our theorem for $r \ne 2$.

\section{Preliminaries}

{\it Notation.}
 
Let $M$ be a complex manifold.
For a cohomology class $x \in H^*(M,{\Bbb Z})$,
$[x]_i \in H^{2i}(X,{\Bbb Z})$ denotes the $2i$-th
component of $x$.

Let $p:X \to \Spec({\Bbb C})$ be an abelian surface  or a K3 surface
over ${\Bbb C}$.
We denote the projection $S \times X \to S$ by $p_S$.
In this paper, we identify a divisor class $D$ with associated 
line bundle ${\cal O}_X(D)$.

\subsection{Mukai lattice}

We shall recall the Mukai lattice [Mu3].
\begin{defn}
We define a symmetric bilinear form on 
$H^{ev}(X,{\Bbb Z}):=\oplus _i H^{2i}(X,{\Bbb Z})$:
\begin{align*}
\langle x, y \rangle:=&-\int_X(x^{\vee} y)\\
=&\; \int_X(x_1y_1-x_0y_2-x_2y_0)
\end{align*}
where $x=x_0+x_1+x_2, x_i \in  H^{2i}(X,{\Bbb Z})$
(resp. $y=y_0+y_1+y_2, y_i \in  H^{2i}(X,{\Bbb Z})$) and
$\vee:H^{ev}(X,{\Bbb Z}) \to H^{ev}(X,{\Bbb Z})$
be the homomorphism sending $x$
to $x_0-x_1+x_2 \in H^{ev}(X,{\Bbb Z})$.
We call this lattice Mukai lattice.
\end{defn}
For a coherent sheaf $E$ on $X$,
\begin{align*}
v(E):=& \ch(E)\sqrt{\td_X}\\
=& \ch(E)(1+\varepsilon \omega)
\end{align*}
is the Mukai vector of $E$, 
where $\omega$ is the fundamental class of $X$ and 
$\varepsilon=0,1$ according as $X$ is of type abelian or K3.
Then Riemann-Roch theorem is written as follows:
\begin{equation}\label{eq:RR}
\chi(E,F)=-\langle v(E),v(F) \rangle,
\end{equation}
where $E$ and $F$ are coherent sheaves on $X$.

For an abelian surface $X$,
Mukai lattice has a decomposition 
\begin{equation}\label{eq:decomp}
\begin{split}
H^{ev}(X,{\Bbb Z}) & =H^2(X,{\Bbb Z}) \oplus H^0(X,{\Bbb Z}) \oplus
H^4(X,{\Bbb Z})\\
&=U^{\oplus 4},
\end{split}
\end{equation}
where $U$ is the hypabolic lattice.

\subsection{Moduli of stable sheaves}
Let $X$ be an abelian or a K3 surface, and    
$H$ an ample line bundles on $X$.
For $v \in H^{ev}(X,{\Bbb Z})$, 
let $M_H(v)$ be the moduli of stable sheaves
of Mukai vector $v$.
By Mukai [Mu2], $M_H(v)$ is smooth of dimension $
\langle v^2 \rangle+2$ and has a symplectic structure.
We set
$$
v^{\perp}:=\{x \in H^{ev}(X,{\Bbb Z})| \langle v, x \rangle=0 \}.
$$

Let $\theta_v: v^{\perp} \to H^2(M_H(v),{\Bbb Z})$
be the homomorphism such that
\begin{equation}
\theta_v(x):=-\frac{1}{\rho}\left[p_{M_H(v)*}((\ch {\cal E})\sqrt{\td_X}
x^{\vee})\right]_1
\end{equation}
where ${\cal E}$ is a quasi-universal family of similitude 
$\rho$.
For a line bundle $L$ on $X$,
let $T_L:H^{ev}(X,{\Bbb Z}) \to H^{ev}(X,{\Bbb Z})$
be the homomorphism sending $x$ to $x\ch(L)$.
Then $T_L$ is an isometry of Mukai lattice and satisfies that 
\begin{equation}\label{eq:isomT}
\theta_{T_L(v)}(T_L(x))=\theta_v(x), 
\end{equation}
for $x \in v^{\perp}$.

From now on, we assume that $X$ is an abelian surface.
Let $\widehat{X}$ be the dual abelian variety of $X$ and
${\cal P}$ the Poincar\'{e} line bundle on
$\widehat{X} \times X$.
For an element $E_0 \in M_H(v)$,
let $\alpha:M_H(v) \to X$ be the morphism
sending $E \in M_H(v)$ to 
$\det p_{\widehat{X}!}((E-E_0)\otimes({\cal P}-{\cal O}_{\widehat{X} \times X}))
\in \Pic^0(\widehat{X})=X$,
and $\det:M_H(v) \to \widehat{X}$
the morphism sending $E$ to $\det E \otimes \det E_0^{\vee}
\in \widehat{X}$.
We set ${\frak a}:=\alpha \times \det$.
Then the following hold [Y2, Thm. 3.1, 3.6].

\begin{thm}\label{thm:H2}
Let $v=r+\xi+a \omega$, $\xi \in H^2(X,{\Bbb Z})$
be a Mukai vector such that $r>0$ and $r+\xi$ is primitive.
We assume that $\dim M_H(v)=\langle v^2 \rangle+2 \geq 6$.
Then for a general ample line bundle $H$, the following holds.\newline
$(1)$ $\theta_v$ is injective.
\newline
$(2)$
${\frak a}$ is the albanese map.
\newline
$(3)$
\begin{equation}\label{eq:H2}
\begin{split}
H^2(M_H(v),{\Bbb Z}) &=
 \theta_v(v^{\perp})
\oplus {\frak a}^* H^2(X \times \widehat{X},{\Bbb Z}).
\end{split}
\end{equation}
\end{thm}

Let $v$ be the Mukai vector in Theorem \ref{thm:H2}.
We set $K_H(v):={\frak a}^{-1}((0,0))$.
We shall construct an \'{e}tale covering 
such that ${\frak a}$ becomes trivial. 
By Theorem \ref{thm:period},
it will become a Bogomolov decomposition of $M_H(v)$.

Let ${\mathbf D}(X)$ and ${\mathbf D}(\widehat{X})$ 
be the derived categories of
$X$ and $\widehat{X}$ respectively.
Let ${\cal S}:{\mathbf D}(X) \to {\mathbf D}(\widehat{X})$ be the Fourier-Mukai
transform in [Mu4], that is,
${\cal S}(F):={\mathbf R}p_{\widehat{X}*}({\cal P} \otimes F),
F \in {\mathbf D}(X)$.
Then $\alpha(E)=\det {\cal S}(E) \otimes (\det{\cal S}(E_0))^{-1}$.
For a line bundle $L$ on $X$,
we set $\hat{L}:=\det({\cal S}(L))$.
Then the following relations hold.
\begin{lem}
\begin{equation}
\begin{split}
& \phi_{\hat{L}} \circ \phi_L=-\chi(L)1_X,\\
&  \phi_{L} \circ \phi_{\hat{L}}=-\chi(L)1_{\widehat{X}}. 
\end{split}
\end{equation}
\end{lem}
 
\begin{pf}
By [Mu4, Prop. 1.21],
$c_1(\hat{\hat{L}})=c_1((-1)^*L)=c_1(L)$ and
$(c_1(\hat{L})^2)=(c_1(L)^2)$.
So it is sufficient to prove the first equality.
By [Mu1, (3.1)], we see that 
${\cal S}(L) \otimes {\cal P}_x={\cal S}(T^*_{-x}L)
={\cal S}(L \otimes {\cal P}_{\phi_L(-x)})=
T^*_{\phi_L(-x)}({\cal S}(L))$.
Hence we get that 
$\widehat{L} \otimes {\cal P}_{\chi(L)x}=
T_{\phi_L(-x)}^*(\widehat{L})=
\widehat{L} \otimes {\cal P}_{\phi_{\widehat{L}}\circ \phi_L(-x)}$.
Therefore the first equality holds.
\end{pf}
We define a morphism
$\Phi:K_H(v) \times X \times \widehat{X} \to M_H(v)$
by
$\Phi(E,x,y):=T_x^*(E) \otimes {\cal P}_y$.

\begin{lem}\label{lem:alb}
Let $L$ be a line bundle on $X$ such that $c_1(L)=c_1$. Then,
\begin{equation}
\begin{split}
\alpha(T_x^*(E) \otimes {\cal P}_y) &=-ax+\phi_{\hat{L}}(y)\\
\det(T_x^*(E) \otimes {\cal P}_y) &=\phi_L(x)+ry.
\end{split}
\end{equation}
\end{lem}
\begin{pf}
We shall only prove the first equality.
By [Mu1, (3.1)],
we see that 
${\cal S}(T_x^*(E) \otimes {\cal P}_y)=T_y^*({\cal S}(T_x^*(E)))
=T_y^*({\cal S}(E) \otimes {\cal P}_{-x})=
T_y^*({\cal S}(E)) \otimes{\cal P}_{-x}$.
Hence $\det({\cal S}(T_x^*(E) \otimes {\cal P}_y))=
T_y^*(\det {\cal S}(E))\otimes{\cal P}_{-\chi(E)x}$. 
Since $c_1( {\cal S}(E))=c_1( {\cal S}(L))$,
$\alpha(T_x^*(E) \otimes {\cal P}_y)=
\phi_{\det {\cal S}(E)}(y)-\chi(E)x=
\phi_{\hat{L}}(y)-\chi(E)x$.
By \eqref{eq:RR}, $\chi(E)=-\langle v({\cal O}_X),v(E) \rangle=a$,
and hence we get the first equality.
\end{pf}

Let $\tau:X \times \widehat{X} \to X \times \widehat{X}$
be a homomorphism sending
$(x,y)$ to
$(rx-\phi_{\hat{L}}(y),-\phi_L(x)-ay)$.
By Lemma \ref{lem:alb},
${\frak a} \circ \Phi \circ (1_{K_H(v)} \times \tau)
(E,x,y)=(nx,ny)$, where $n=\langle v^2 \rangle/2$.
Let $\nu:X \times \widehat{X} \to X \times \widehat{X}$
be the $n$ times map and
we shall consider the fiber product
\begin{equation}
\begin{CD}
M_H(v) \times_{X \times \widehat{X}}X \times \widehat{X} @>>> M_H(v)\\
@VVV @VV{{\frak a}}V \\
X \times \widehat{X} @>{\nu}>> X \times \widehat{X}
\end{CD}
\end{equation}
Then $\Phi \circ (1_{K_H(v)} \times \tau)$ and the projection
$K_H(v) \times X \times \widehat{X} \to  X \times \widehat{X}$
defines a morphism to the fiber product.
We can easily show that this morphism is injective,
and hence it is an isomorphism.
 
\begin{rem}
If $(c_1^2)/2$ and $r$ are relatively prime,
then [Y2, Prop. 4.1] implies that 
$M_H(v) \cong \widehat{X} \times \det^{-1}(0)$.
We shall consider the pull-back of 
${\frak a}: M_H(v) \to X \times \widehat{X}$ by the morphism
sending $(x,y)$ to $(nx,y)$.
Then we get $M_H(v) \times_{ X \times \widehat{X}} X \times \widehat{X}
\cong K_H(v) \times  X \times \widehat{X}$.
\end{rem}

For simplicity,
we also denote the homomorphism
$v^{\perp} \to H^2(M_H(v),{\Bbb Z}) \to H^2(K_H(v),{\Bbb Z})$
by $\theta_v$:
\begin{equation}
\theta_v(x)=-\frac{1}{\rho}\left[p_{K_H(v)*}(\ch({\cal E}_{|K_H(v) \times X})
x^{\vee})\right]_1.
\end{equation}

\subsection{Beauville's bilinear form}
Let $M$ be an irreducible symplectic manifold
of dimension $n$.
Beauville [B] constructed a primitive symmetric bilinear
form 
\begin{equation}
B_M:H^2(M,{\Bbb Z}) \times H^2(M,{\Bbb Z}) \to {\Bbb Z}.
\end{equation}
Up to multiplication by positive constants,
$q_M(x):=B_M(x,x)$ satisfies that
\begin{equation}
q_M(x)=\frac{n}{2}\int_M \phi^{n-1} \overline{\phi}^{n-1} x^2
+(1-n)\int_M \phi^{n} \overline{\phi}^{n-1} x
\int_M \phi^{n-1}  \overline{\phi}^{n} x,
\end{equation}
where $\phi$ is a holomorphic 2 form with 
$\int_M \phi^{n}\overline{\phi}^{n}=1$.
For $\lambda, x \in H^2(M,{\Bbb C})$,
the following relation holds [B, Thm. 5].
\begin{equation}\label{eq:Bform}
v(\lambda)^2q_M(x)=q_M(\lambda)\left[
(2n-1)v(\lambda)\int_M \lambda^{2n-2} x^2
-(2n-2)\left(\int_M \lambda^{2n-1} x\right)^2 \right],
\end{equation}
where $v(\lambda)=\int_M \lambda^{2n}$.

\section{Proof of Theorem \ref{thm:period}}

\subsection{Generalized Kummer variety}
 
In this subsection, we shall recall Beauville's results [B]
on generalized Kummer varieties.
Then Theorem \ref{thm:period} for $r=1$ follows from his results and
simple calculations.
Let $X$ be an abelian surface.
Let $\pi:X^n \to X^{(n)}$
be the $n$-th symmetric product of $X$.
We set $X^{[n]}:=\Hilb_X^n$.
Let $\gamma:X^{[n]} \to X^{(n)}$ be the Hilbert-Chow morphism.
Let $\sigma:X^{(n)} \to X$ be the morphism sending 
$(x_1,x_2,\dots,x_n) \in X^{(n)}$ to $\sum_{i=1}^n x_i \in X$.
Then ${\frak a}:X^{[n]} \to X^{(n)} \to X$ is the albanese map of 
$X^{[n]}$.
If $n=2$, then ${\frak a}^{-1}(0)$ is the
Kummer surface associated to $X$ 
and if $n \geq 3$, then
$K_{n-1}:={\frak a}^{-1}(0)$ is the generalized Kummer variety
constructed by Beauville [B].
$K_{n-1}$ is an irreducible symplectic manifold of
dimension $2(n-1)$.

We assume that $n \geq 3$.
For integers $i,j,k$, we set
$\Delta^{i,j}:=\{(x_1,x_2,\dots,x_n) \in X^n| x_i=x_j \}$,
$\Delta^{i,j,k}:=\Delta^{i,j} \cap \Delta^{j,k}$.
We set $X_*^n:=X^n \setminus \cup_{i<j<k} \Delta^{i,j,k}$,
$X_*^{[n]}:=X^{[n]} \setminus_{i<j<k}\gamma^{-1}(\pi(\Delta^{i,j,k}))$.
We set $N:=\{(x_1,x_2,\dots,x_n)| x_1+x_2+\dots+x_n=0 \}$,
$N_*:=N \cap X_*^n$, $(K_{n-1})_*:=K_{n-1} \cap X_*^{[n]}$,
$\delta^{i,j}:=\Delta^{i,j} \cap N$.
Since $n \geq 3$,
$\delta^{i,j}$ is connected,
indeed, it is isomorphic to $X^{n-2}$.
Let $\beta:B_{\Delta}(X_*^n) \to X_*^n$ be the blow-up of
$X_*^n$ along $\delta:=\cup_{i<j} \Delta^{i,j}$
and set $B_{\delta}(N_*)=\beta^{-1}(N_*)$.
Let $E^{i,j}:=\beta^{-1}(\Delta^{i,j})$ 
be the exceptional divisor of $\beta$
and $e^{i,j}:=\beta^{-1}(\Delta^{i,j} \cap N)$.

\begin{equation}
\begin{CD}
B_{\delta}(N_*) @>>> B_{\Delta}(X_*^{n}) @>{\beta}>> X^n\\
@VV{\pi'}V @VV{\pi'}V @VV{\pi}V \\
(K_{n-1})_* @>>> X_*^{[n]} @>>{\gamma}> X^{(n)}\\
@.@.@VV{\sigma}V\\
@.@.X
\end{CD}
\end{equation}
We shall first describe $H^2(K_{n-1},{\Bbb Z})$.

\begin{lem}
$$
H^2(K_{n-1},{\Bbb Z}) \cong H^2(X,{\Bbb Z}) \oplus {\Bbb Z}e.
$$
\end{lem}
\begin{pf}
By [B, Prop. 8],
$H^2(K_{n-1},{\Bbb Q}) \cong H^2(X,{\Bbb Q}) \oplus {\Bbb Q}\,e$.
Since $\varphi:H^2(X,{\Bbb Z}) \to H^2(K_{n-1},{\Bbb Z})
\to H^2(B_{\delta}(N_*),{\Bbb Z})$ is injective and
$\im \varphi \subset \beta^*(H^2(N,{\Bbb Z}))$,
we shall prove that the image of 
$f:H^2(X,{\Bbb Z}) \to H^2(N,{\Bbb Z})^{{\frak S}_n}$ is a
primitive submodule of
$H^2(N,{\Bbb Z})$.
Let $\phi:X \times X \to N$ be the morphism such that 
$\phi((x,y))=(x,y,0,\dots,0,-x-y) \in N$.
We shall consider the composition $g:H^2(X,{\Bbb Z}) \to 
H^2(X \times X,{\Bbb Z})$.
Let $\alpha_i \wedge \alpha_j,i<j$ be the basis of 
$H^2(X,{\Bbb Z})=\wedge^2 H^1(X,{\Bbb Z})$.
Then we see that
$g^*(\alpha_i \wedge \alpha_j)=2 p_1^*(\alpha_i \wedge \alpha_j)+
2p_2^*(\alpha_i \wedge \alpha_j)+(p_1^*\alpha_i \wedge p_2^*\alpha_j
-p_1^* \alpha_j \wedge p_2^*\alpha_i)$.
Hence $\im g$ is a primitive subspace of $H^2(X \times X,{\Bbb Z})$.
Therefore $\im f$ is primitive.
\end{pf}

We shall prove that $\theta_v$ preserves the bilinear forms.
Since $v=1-n \omega$,
we get that
$v^{\perp}=H^2(X,{\Bbb Z})\oplus {\Bbb Z}(1+n \omega)$. 
For $\alpha=x+k(1+n \omega),
x \in H^2(X,{\Bbb Z})$,
simple calculations show that
\begin{equation}
\begin{split}
\langle \alpha^2 \rangle &=(x^2)-k^2(2n),\\
\theta_v(\alpha)&=\sum_i p_i^* (x)+k e.
\end{split}
\end{equation}
Hence we shall prove that
\begin{equation}\label{eq:q}
q_{K_{n-1}}(\theta_v(\alpha))=(x^2)-k^2 (2n).
\end{equation}

We shall choose $l \in H^2(X,{\Bbb Z})$ with $(l^2) \ne 0$.
We set $\omega=\alpha_1 \wedge \alpha_2 \wedge \alpha_3 \wedge \alpha_4$.
Then the $\oplus_{i<j}p^*_i H^{ev}(X,{\Bbb Z}) \otimes
p_j^* H^{ev}(X,{\Bbb Z})$-component of
$\sigma^* (\omega)=\sum_{i,j,k,m}p_i^* \alpha_1
\wedge p_j^* \alpha_2 \wedge p_k^* \alpha_3 \wedge p_m^*\alpha_4$
is $\sum_i p_i^*(\alpha_1 \wedge \alpha_2 \wedge \alpha_3 \wedge \alpha_4)
+\sum_{i \ne j}p_i^*(\alpha_1 \wedge \alpha_2)
 \wedge p_j^*(\alpha_3 \wedge \alpha_4)$.
Since the $\oplus_{i<j}p^*_i H^{ev}(X,{\Bbb Z}) \otimes
p_j^* H^{ev}(X,{\Bbb Z})$-component of 
$\sum_{i<j}(\Delta^{i,j}
-p_i^* \omega -p_j^* \omega)$ is $\sum_{i \ne j}p_i^*(\alpha_1 \wedge \alpha_2)
 \wedge p_j^*(\alpha_3 \wedge \alpha_4)$,
those of $\mu:=\sum_i p_i^* \omega
+\sum_{i<j}(\Delta^{i,j}
-p_i^* \omega -p_j^* \omega)$ and $\sigma^* (\omega)$ are the same.
Hence we see that
\begin{equation}
\begin{split}
\int_{K_{n-1}}(\theta_v(l))^{2n-4}(\theta_v(x))^2 &=
\frac{1}{n!} \int_N (p_1^* l+\dots+ p_n^*l)^{2n-4}
(p_1^*x+\dots +p_n^* x)^2\\
&=\frac{1}{n!} \int_{X^n} (p_1^* l+\dots +p_n^*l)^{2n-4}
(p_1^*x+\dots +p_n^* x)^2 \mu\\
&=\frac{1}{n!}\left\{\frac{n(n-1)}{2}
 \int_{X^{n-1}} (2p_1^* l+p_2^*l+\dots +p_n^*l)^{2n-4}
(2p_1^*x+p_2^*x+\dots +p_n^* x)^2 \right.\\
& \quad \left.-n(n-2)\int_{X^{n-1}} (p_1^* l+p_2^*l+\dots +p_n^*l)^{2n-4}
(p_1^*x+p_2^*x+\dots +p_n^* x)^2
\right\}\\
&=\frac{(2n-2)!n^2}{n!2^{n-1}}\left(
\frac{1}{2n-3}(l^2)^{n-2}(x^2)+\frac{2n-4}{2n-3}(l^2)^{n-2}(l,x)^2 
\right).
\end{split}
\end{equation}
In the same way, we see that
\begin{equation}
\begin{split}
\int_{K_{n-1}}(\theta_v(l))^{2n-3}(\theta_v(x)) &=
\frac{(2n-2)!n^2}{2^{n-1}n!}(l^2)^{n-1}(l,x),\\
\int_{K_{n-1}}(\theta_v(l))^{2n-2}&=
\frac{(2n-2)!n^2}{2^{n-1}n!}(l^2)^{n}.
\end{split}
\end{equation}
By \eqref{eq:Bform}, we obtain that
\begin{equation}\label{eq:D}
q_{K_{n-1}}(\theta_v(x))=\frac{(x^2)}{(l^2)}q_{K_{n-1}}(\theta_v(l)).
\end{equation}

We shall next compute $q_{K_{n-1}}(e)$.
We note that $\gamma^{-1}(\pi(\Delta^{1,2}))=2e$.
Let $\iota:H^2(X,{\Bbb Z}) \to H^2(X^{(n)},{\Bbb Z})$ be the homomorphism 
such that $\pi^*(\iota(x))=\sum_i p_i^* x \in H^2(X^n,{\Bbb Z})$. 
Since the Poincar\'{e} dual of
$(\iota(l))^{2n-4}\iota(x)$, $x \in H^2(X,{\Bbb Z})$
(resp. $(\iota(l))^{2n-4}$) is
a 2 cycle (resp. 4 cycle) of $X^{(n)}$, 
the intersection with $\pi(\Delta_{1,2})$ is $0$
(resp. finite points).

Let $\sigma':X^{n-1}=\Delta_{1,n} \to X$ be the restriction
of $\sigma$ to the diagonal $\Delta_{1,n}$. 
We set
$\mu'=2^4 p_1^*(\omega)+\sum_{i=2}^{n-1}p_i^*(\omega)+
\sum_{1<i<j \leq n-1}
(\Delta^{i,j}-p_i^*(\omega)-p_j^*(\omega))
+4\sum_{i=2}^{n-1}(\Delta_{1,i}-p^*_1 \omega-p^*_i\omega)
$.
In the same way, we see that the 
 $\oplus_{i<j}p^*_i H^{ev}(X,{\Bbb Z}) \otimes
p_j^* H^{ev}(X,{\Bbb Z})$-components of
$\sigma^* (\omega)$ and $\mu'$ 
are the same.
Since $E_{i,j}$ is the exceptional divisor
of $\beta$,
$\beta_*((E_{i,j})_{|E_{i,j}})=-\Delta^{i,j}$.
Hence we see that 

\begin{equation}
\begin{split}
\int_{K_{n-1}}(\theta_v(l))^{2n-4}\theta_v(x) e=& 0,\\
\int_{K_{n-1}}(\theta_v(l))^{2n-4} e^2 =&
-\frac{1}{n!} \int_N
(p_1^* l+\dots+ p_n^*l)^{2n-4}
(\sum_{i<j} \Delta^{i,j})\\
=&-\frac{n(n-1)}{2 n!} \int_{X^{n-1}} 
(2p_1^* l+p_2^*l+\dots +p_{n-1}^*l)^{2n-4}
\mu'\\
=&\frac{(2n-2)! n^2}{2^{n-1}n!}\frac{(-2n)}{2n-3}.
\end{split}
\end{equation}
Thus $e$ is orthogonal to $H^2(X,{\Bbb Z})$ and
\begin{equation}\label{eq:e}
q_{K_{n-1}}(e)=\frac{-2n}{(l^2)}q_{K_{n-1}}(\theta_v(l)).
\end{equation}
By \eqref{eq:D} and \eqref{eq:e},
we get \eqref{eq:q}. 

\begin{prop}[Beauville]
For $v=1-n \omega$, $n \geq 2$,
\begin{equation}
\theta_v:v^{\perp} \to H^2(K_H(v),{\Bbb Z})
\end{equation}
is an isometry of Hodge structures.
\end{prop}

\subsection{General cases}

We shall first treat the case
where $X$ is a product of two elliptic curves,
and by deformation arguments,
we shall treat general cases.
Let $X$ be an abelian surface which is a product of
two elliptic curves $C_1, C_2$. 
Let $f_i$, $i=1,2$ be the ample generator of
$H^2(C_i,{\Bbb Z})$.
Let $(r,d)$ and $(r_1,d_1)$ be pairs of integers
such that $r>r_1>0$ and
$d r_1-r d_1=1$.
We set $v=r+d f_2-(r-r_1)nf_1-(d-d_1)n \omega$.
We shall choose an ample divisor $H=f_2+m f_1, m \gg 0$. 
In [Y2, sect. 3.2],
we constructed an immersion
$B_{\Delta}(X^n)/{\frak S}_n \to M_H(v)$ which is an isomorphism in codimension 1.
Hence we get an isomorphism  
$H^2(M_H(v),{\Bbb Z}) \to H^2(B_{\Delta}(X^n),{\Bbb Z})^{{\frak S}_n}$.
It also induces a birational map
$(K_{n-1})_* \to K_H(v)$ which is an isomorphism in codimension 1.
Thus we get an isomorphism
$H^2(K_H(v),{\Bbb Z}) \to H^2(K_{n-1},{\Bbb Z})$. 
By the injective homomorphism
$H^2(K_H(v),{\Bbb Z}) \to H^2(B_{\delta}(N_*),{\Bbb Z})$,
we shall regard $H^2(K_H(v),{\Bbb Z})$ as an submodule of
$H^2(B_{\delta}(N_*),{\Bbb Z})$.

We set $x=x_1+x_2 f_1 +x_3 f_2+x_4 \omega+D$,
$D \in H^1(C_1,{\Bbb Z}) \otimes H^1(C_2,{\Bbb Z})$.
We assume that $r \geq 2r_1$.
In the notation of [Y2, sect. 3.1],
$\theta_v(x)=-\kappa_2(x^{\vee})$.
Hence [Y2, (3.19), (3,20)] are written down
as follows:
\begin{equation}
\theta_v(x)
=y_1(\sum_{i=1}^n p_i^*f_2) +y_2 (\sum_{i=1}^n p_i^* f_1)
+y_3(\sum_{i<j} E^{i,j})+\sum_{i=1}^n p_i^* D,
\end{equation}
where
\begin{equation}
\begin{cases}
y_1=dx_1-r x_3\\
y_2=-(d-d_1)x_2+(r-r_1) x_4-n((d-2d_1)x_1-(r-2r_1)x_3)\\
y_3=-d_1x_1+r_1x_3\\
y_4=dx_2-rx_4+n((d-d_1)x_1-(r-r_1)x_3).
\end{cases}
\end{equation}
By simple calculations, we get that
\begin{equation}
\begin{cases}
x_1=r_1 y_1+r y_3\\
x_2=-nr_1 y_1-ry_2-n(r+r_1)y_3-(r-r_1)y_4\\
x_3=d_1y_1+d y_3\\
x_4=-n d_1y_1-d y_2-n(d+d_1)y_1-d_2 y_4.
\end{cases}
\end{equation}
By the definition of $v^{\perp}$,
$x$ belongs to $v^{\perp}$
if and only if $y_4=0$.
Hence we obtain that
\begin{equation}
\begin{split}
\langle x^2 \rangle &=
2x_2x_3-2x_1x_4+(D^2)\\
&=2 y_1 y_2+(D^2)-2ny_3^2\\
&=q_{K_H(v)}(\theta_v(x)).
\end{split}
\end{equation}

\begin{lem}\label{lem:isom2}
Under the same assumptions on $X$,
let $v=r+(d f_2 +s f_1)+a \omega$
be a Mukai vector such that
$(r,d)=1$ and $\langle v^2 \rangle=2n \geq 4$.
We set $H:=f_2+m f_1$, where $m$ is a sufficiently
large integer. 
Then
$K_H(v)$ is an irreducible symplectic manifold and
\begin{equation}
\theta_v:v^{\perp} \to H^2(K_H(v),{\Bbb Z})
\end{equation}
is an isometry of Hodge structures for $n \geq 3$.
\end{lem}
\begin{pf}
We shall choose a pair of integers $(r_1,d_1)$
such that $r>r_1>0$ and $d r_1-rd_1=1$.
We first assume that $r \geq 2 r_1$.
Since $\langle v^2 \rangle=2(ds-ra)=2n$,
there are integers $s_1,a_1$ such that
\begin{equation}
\begin{cases} 
s=n r_1+s_1 r,\\
a=n d_1+a_1 d.
\end{cases}
\end{equation}
We set $v':=v \ch({\cal O}_X(-(n+s_1)f_1))$.
Then we see that
$v'=r+d f_2-(r-r_1)n f_1-(d-d_1)n \omega$.
Since $K_H(v')$ is an irreducible symplectic manifold and 
$\theta_{v'}$ is an isometry,
combining \eqref{eq:isomT}, the assertions hold for $K_H(v)$ with
$r \geq 2r_1$.
If $r<2r_1$, then we shall replace $v$ by $v^{\vee}$.
Since $\theta_{v^{\vee}}(x^{\vee})=-\theta_v(x)$, 
$x \in H^{ev}(X,{\Bbb Z})$, this case can be reduced to 
the first case.
\end{pf}

We shall treat general cases.
Twisting by some ample line bundles,
we may assume that $\xi$ belongs to the ample cone. 
The following argument is the same as that in [Y2, Prop. 3.3].
Let $f:({\cal X},{\cal L}) \to T$ be a family of polarized abelian surfaces
over a connected curve $T$ such that $\NS({\cal X}_t)={\Bbb Z}{\cal L}_t$
for some $t \in T$.
Let $v:=r+d {\cal L}+a \omega \in R^{ev}f_*{\Bbb Z}$
be a family of Mukai vector such that $(r,d)=1$.
By [Y2, Prop. 3.3], we can construct a proper and smooth family
of moduli spaces ${\cal M}_{{\cal X}/T}(v) \to T$ and a family of
albanese maps ${\frak a}_T:{\cal M}_{{\cal X}/T}(v) \to {\cal X}\times_T
Pic_{{\cal X}/T}^0$.
Since ${\cal X} \to T$ is projective,
we can also construct a family of homomorphisms
$(\theta_v)_t:(v^{\perp})_t \to H^2(({\cal M}_{{\cal X}/T}(v))_t,{\Bbb Z})$.
Let $0_T:T \to {\cal X}\times_T
Pic_{{\cal X}/T}^0$ be the $0$-section of $f$, and we set
${\cal K}_{{\cal X}/T}(v):={\frak a}_T^{-1}(0_T)$.
We assume that $({\cal K}_{{\cal X}/T}(v))_{t_0}, t_0 \in T$ is irreducible
symplectic and $(\theta(v))_{t_0}:(v^{\perp})_{t_0} \to 
H^2(({\cal K}_{{\cal X}/T}(v))_{t_0},{\Bbb Z})$ is an isometry of Hodge structures.
Then every fiber of ${\frak a}$ is irreducible symplectic and
$(\theta_v)_t$ is an isometry of Hodge structures. 
We note that moduli of $(1,n)$-polarized abelian surfaces is irreducible.
Applying this assertion, we can reduce to the case where
$X$ is a product of elliptic curves.
Applying Lemma \ref{lem:isom2},  
we get Theorem \ref{thm:period}.
\qed
\hfill \hfill
\begin{cor}
We set $(v^{\perp})_{alg}:=v^{\perp} \cap
(H^0(X,{\Bbb Z}) \oplus \NS(X)\oplus
H^4(X,{\Bbb Z}))$.
Then $\theta_v$ induces an isometry 
$$
(v^{\perp})_{alg} \to \NS(K_H(v)).
$$
\end{cor}

The following example is similar to [Mu5, 5.17]. 

\begin{ex}\label{ex:ex1}
Let $X$ be an abelian surface with $\NS(X)={\Bbb Z}H$,
$(H^2)=2$.
We set $v=2+H-2 \omega$.
Then $M_H(v)$ is a variety of dimension $12$.
It is easy to see that
$v^{\perp}$ is generated by
$\alpha:=1+\omega$ and $\beta=H+\omega$.
Since
\begin{alignat}3
\langle \alpha^2 \rangle=-2, \quad & \langle \alpha,\beta \rangle=-1, 
\quad & \langle \beta^2 \rangle=2,
\end{alignat} 
$\NS(K_H(v))$ is indecomposable.
Hence $M_H(v)$ is not birationally equivalent to
$\widehat{Y} \times \Hilb_Y^{5}$ for any $Y$.
\end{ex}

\section{The case of $\langle v^2 \rangle=4$}

In this section, we shall treat the remaining case,
that is, $\langle v^2 \rangle =4$.
In this case, $K_H(v)$ is a K3 surface.
We shall determine this K3 surface.
Let $v=r+\xi+a \omega$,
$\xi \in H^2(X,{\Bbb Z})$ be a Mukai vector such that
$r+\xi$ is primitive and $\langle v^2 \rangle =4$.
Replacing $v$ by $v \ch(H^{\otimes m})$, $m \gg 0$,
we may assume that $\xi$ belongs to the ample cone.
Let $\iota:X \to X$ be the $(-1)$-involution of $X$
and $x_1,x_2,\dots,x_{16}$ the fixed points of $\iota$.
Let $\pi:\widetilde{X} \to X$ be the blow-ups of $X$ at 
$x_1,x_2,\dots,x_{16}$ and $E_1,E_2,\dots,E_{16}$
the exceptional divisors of $\pi$. 
Let $q_1:X \to X/\iota$ be the quotient map.
Then the morphism $q_1 \circ \pi:\widetilde{X} \to X/\iota$ factors through
the quotient $\widetilde{X}/\iota$ of $\widetilde{X}$ by $\iota$ 
: $\widetilde{X} \overset{q_2} \to \widetilde{X}/\iota
\overset{\varpi} \to X/\iota$.
$\Km(X):=\widetilde{X}/\iota$ is the Kummer surface associated to $X$
and $\varpi:\widetilde{X}/\iota \to X/\iota$ is the minimal resolution of 
$X/\iota$.
We set $C_i:=q_2(E_i)$, $i=1,2,\dots,16$.
 
We may assume that $H$ is symmetric,
that is, $\iota^*H=H$.
Then $H$ has a $\iota$-linearization.
Hence $H^{\otimes 2}$ descent to an ample line bundle $L$ on
$X/\iota$.
Then $L_m:=\varpi^*(L^{\otimes m})(-\sum_{i=1}^{16}C_i)$,
$m \gg0$ is an ample line bundle on $\Km(X)$.
We shall fix a sufficiently large integer $m$.
Let $w=r+c_1+b \omega \in H^{ev}(\Km(X),{\Bbb Z})$
be an isotropic Mukai vector.
We shall consider moduli space $M_{L_m}(w)$.
By Mukai [Mu3], $M_{L_m}(w)$ is not empty.
By our assumption on $L_m$, $M_{L_m}(w)$ consists of $\mu$-stable 
sheaves.
Indeed, for a $\mu$-semi-stable vector bundle $F$ of
$v(F)=w+k \omega, k \geq 0$,
$q^*_2(F)$ is a $\mu$-semi-stable vector bundle on $X$ with
respect to $\pi^*(H^{\otimes 2m})(-2\sum_{i=1}^{16}E_i)$.
Since $m$ is sufficiently large and $H$ is a general
ample line bundle,
$q_2^*(F)$ is a $\mu$-stable vector bundle.
Hence $F$ is a $\mu$-stable vector bundle,
which implies that $M_{L_m}(w)$ consists of $\mu$-stable sheaves.
Since $\dim M_{L_m}(w)=\langle w^2 \rangle+2=2$,
every member of $M_{L_m}(w)$ is locally free.
Moreover general members $F$ of $M_{L_m}(w)$ are 
rigid on each $(-2)$-curves $C_i$.
\begin{lem}
We set
$N(w,i):=\{F \in M_{L_m}(w)|\text{$F_{|C_i}$ is not rigid}\}$.
Then $N(w,i)$ is not empty if and only if $r|\deg(F_{|C_i})$.
Moreover if $N(w,i)$ is not empty, then $N(w,i)$ is a rational curve.
\end{lem}
\begin{pf}
We assume that $F_{|C_i}$ is not rigid.
We set $F_{|C_i}=\oplus_{j=1}^k{\cal O}_{C_i}(a_j)^{\oplus n_j}$,
$a_1<a_2<\dots<a_k$.
Let $F':=\ker (F \to {\cal O}_{C_i}(a_1)^{\oplus n_1})$
be the elementary transformation of $F$ along
${\cal O}_{C_i}(a_1)^{\oplus n_1}$.
Then $v(F')=v(F)-n_1(C_i-(a_1+1)\omega)$.
Hence we see that
$\langle v(F')^2\rangle=-2 n_1(\sum_{j \geq 2}n_j(a_j-a_1-1))$.
By the choice of $L_m$,
$F'$ is also $\mu$-stable.
Hence $-2 \leq -2 n_1(\sum_{j \geq 2}n_j(a_j-a_1-1))$.
Since $F$ is not rigid,
$\sum_{j \geq 2}n_j(a_j-a_1-1) > 0$.
Thus $n_1=\sum_{j\geq 2} n_j(a_j-a_1-1)=1$.
Therefore we get that
$F_{|C_i} \cong {\cal O}_{C_i}(a_1)
\oplus {\cal O}_{C_i}(a_1+1)^{\oplus (r-2)} \oplus  {\cal O}_{C_i}(a_1+2)$.
In this case, $\langle v(F')^2 \rangle =-2$, and hence
$F'$ is a unique stable vector bundle
of $v(F')=v(F)-n_1(C_i-(a_1+1)\omega)$.
It is not difficult to see that the choice of inverse transformations
is parametrized by ${\Bbb P}^1$. Therefore $N(w,i)$ is a rational curve.
\end{pf}

We shall consider the pull-back $q^*_2(F)$ of a general member $F$.
Since $F_{|C_i}$, $1 \leq i \leq 16$ are rigid,
replacing $q_2^*(F)$ by
$q_2^*(F)(\sum_{i=1}^{16} s_i E_i)$,
we may assume that $q^*_2(F)_{|E_i} \cong 
{\cal O}_{E_i}(-1)^{\oplus k_i}  
\oplus {\cal O}_{E_i}^{\oplus (r-k_i)}$.
Let $\phi:q^*_2(F) \to \oplus_{i=1}^{16}{\cal O}_{E_i}(-1)^{\oplus k_i}$
be the quotient map
induced by the quotients $q^*_2(F)_{|E_i} \to {\cal O}_{E_i}(-1)^{\oplus k_i}$.
Then $G:=\ker \phi$ is the elementary transformation of $q_2^*(F)$
along $\oplus_{i=1}^{16}{\cal O}_{E_i}(-1)^{\oplus k_i}$
and $G$ satisfies that $G_{|E_i} \cong {\cal O}_{E_i}^{\oplus r}$.
Hence $\pi_*(G)$ is a stable vector bundle on $X$.
So we get a rational map
$f:M_{L_m}(w) \cdots \to M_H(v)$, where $v=v(\pi_*(G))$.
Since $M_{L_m}(w)$ is a K3 surface, the image of $M_{L_m}(w)$
belongs to a fiber of ${\frak a}$.
Since $q_2^*(F)$ is a stable, and hence a simple vector bundle
and $\iota$ has fixed points,
$\iota$-linearization on $F$ is uniquely determined by
$q_2^*(F)$.
Hence $f$ is generically injective.  
By a simple calculation, we get that 
\begin{equation}\label{eq:v(G)}
\begin{split}
\langle v(G)^2 \rangle &=2r c_2(G)-(r-1)(c_1(G)^2)\\
&=4r c_2(F)-2(r-1)(c_1(F)^2)-\sum_{i=1}^{16}k_i(r-k_i)\\
&=2(\langle w^2 \rangle+2r^2)-\sum_{i=1}^{16}k_i(r-k_i)\\
&=4r^2-\sum_{i=1}^{16}k_i(r-k_i).
\end{split}
\end{equation}
Hence if $\langle v(G)^2 \rangle=4$, then the fiber of ${\frak a}$
is isomorphic to $M_{L_m}(w)$.

Conversely for a Mukai vector $v=r+d N+a \omega \in H^{ev}(X,{\Bbb Z})$
such that $(a)$ $N$ is a $(1,n)$-polarization,
$(b)$ $(r,d)=1$ and $(c)$ $\langle v^2 \rangle =d^2 (N^2)-2r a=4$,
we shall look for such a vector $w \in H^{ev}(\Km(X),{\Bbb Z})$.
We shall divide the problem into two cases.

Case (I). We first assume that 
$r$ is even.
In this case, $d$ must be odd.
By the condition $(c)$, $(N^2)=2n$ is divisible by 4.
Thus $n$ is an even integer.
In this case, replacing $N$ by $N \otimes{\cal P}$ with
${\cal P}^{\otimes 2} \cong {\cal O}_X$,
we may assume that $N$ has a $\iota$-linearization which acts 
trivially on the fibers of $N$ at exactly 4 points (cf. [L-B, Rem. 7.7]).
Replacing the indices, we assume that the 4 points are 
$x_1,x_2,x_3,x_4$.
We set $N_1:=\pi^*(N^{\otimes d})(\frac{r-2}{2}\sum_{i=1}^4 E_i+
\frac{r}{2}\sum_{i \geq 5}E_i)$ and $N_2:=N_1(-r E_1)$.
Then for suitable linearizations,
$N_1$ and $N_2$ descend to line bundles 
$\xi_1$ and $\xi_2$ on $\Km(X)$ respectively.
By simple calculations, we get that
\begin{equation}
\begin{split}
(\xi_1^2) &=d^2 \frac{(N^2)}{2}-2r^2+2r-2\\
&=r(a-2r+2),\\
(\xi_2^2) &= r(a-2r+1).
\end{split}
\end{equation}
We set
\begin{equation}
w:=
\begin{cases}    
 r+\xi_1+\frac{a-2r+2}{2}\omega, \text{ if $a$ is even},\\
r+\xi_2+\frac{a-2r+1}{2}\omega, \text{ if $a$ is odd}.
\end{cases}
\end{equation}
Then we get that $\langle w^2 \rangle =0$.
Let $F$ be a general stable vector bundle of $v(F)=w$.
By the choice of $\xi_1$ and $\xi_2$, 
the restriction of $q_2^*(F)$ or $q_2^*(F)(E_1)$
to $E_i$ is isomorphic to
${\cal O}_{E_i}(-1)^{\oplus k_i} \oplus {\cal O}_{E_i}^{\oplus (r-k_i)}$,
where $k_i=(r-2)/2$ for $1 \leq i \leq 4$ and
$k_i=r/2$ for $i \geq 5$.
Then by \eqref{eq:v(G)},
we get that $\langle v(\pi_*(G))^2 \rangle=4$.
Since $\rk(\pi_*(G))=r$ and $c_1(\pi_*(G))=dN$,
$v(\pi_*(G))$ must be equal to $v$.
Therefore $K_H(v)$ is isomorphic to $M_{L_m}(w)$.
 
Case (II). We assume that $r$ is odd.
Replacing $v$ by $v \ch(N)$, we may assume that $d$ is even.
We set $N_1:=\pi^*(N^{\otimes d})(\frac{r-1}{2}\sum_{i=1}^{16}E_i)$.
Then for a suitable linearization,
$N_1$ descend to a line bundle $\xi$ on $\Km(X)$.
By a simple calculation, we get that 
$(\xi^2)=r(a-2r+4)$.
Since $d$ is even and $r$ is odd, condition $(c)$ implies that
$a$ is an even integer.
We set $v:=r+\xi+\{(a-2r+4)/2\}\omega$.
Then we get that $\langle v^2 \rangle=0$.
In the same way as above, we see that $v(\pi_*(G))=v$,
which implies that $K_H(v) \cong M_{L_m}(w)$.
\begin{thm}
Let $v=r+\xi+a \omega \in H^{ev}(X,{\Bbb Z})$ be a Mukai vector
such that $r>0$, $r+\xi$ is primitive and $\langle v^2 \rangle =4$.
Let $\Km(X)$ be the Kummer surface associated to $X$.
Then there is an isotropic Mukai vector 
$w \in H^{ev}(\Km(X),{\Bbb Z})$ and an ample line bundle $H'$
on $\Km(X)$ such that $K_H(v)$ is isomorphic to $M_{H'}(w)$.
\end{thm} 

\begin{rem}
By the choice of $k_i$,
if $r>2$, then $N(w,i)$ is empty.
Thus $f$ is a morphism.
If $r=2$, then $N(w,i)$, $1 \leq i \leq 4$
is not empty and these closed subset correspond to
the closed subset $N(v,i):=\{G \in K_H(v)|
\text{$G$ is not locally free at $x_i$}\}$.
\end{rem}

\section{appendix}
In this appendix, we shall explain another method to prove 
Theorem \ref{thm:period}.
Let $(X_1, X_2,{\cal P})$ be a triple of surfaces $X_1,X_2$ and 
a coherent sheaf ${\cal P}$ on $X_1 \times X_2$ such that
$K_{X_1}$ and $K_{X_2}$ are trivial, ${\cal P}$ is flat over $X_1$
and $X_2$, and ${\cal P}$ is strongly simple over $X_1$ and $X_2$
(see [Br, sect. 2]).
We denote the projections $X_1 \times X_2 \to X_i$, $i=1,2$ 
 by $p_i$.
Let ${\cal F}_D:{\mathbf D}(X_1) \to {\mathbf D}(X_2)$ be the Fourier-Mukai
transform defined by ${\cal P}$,
that is,
${\cal F}_D(x)={\mathbf R}p_{2*}({\cal P} \otimes p_1^*(x)), 
x \in {\mathbf D}(X_1)$.
Let $\widehat{{\cal F}}_D:{\mathbf D}(X_2) \to {\mathbf D}(X_1)$ 
be the inverse transformation,
that is,
$\widehat{{\cal F}}_D(y)={\mathbf R}\Hom_{p_1}({\cal P},p_2^*(y)),
 y \in {\mathbf D}(X_2)$.
Let ${\cal F}_H:H^{ev}(X_1,{\Bbb Q}) \to H^{ev}(X_2,{\Bbb Q})$ and
$\widehat{{\cal F}}_H:H^{ev}(X_2,{\Bbb Q}) \to H^{ev}(X_1,{\Bbb Q})$ 
be homomorphisms such that
\begin{align}
{\cal F}_H(x)&=p_{2*}((\ch{\cal P})p_1^*\sqrt{\td_{X_1}}p_2^*\sqrt{\td_{X_2}}
p_1^*(x)), x \in H^{ev}(X_1,{\Bbb Q}),\\
\widehat{{\cal F}}_H(y)&=
p_{1*}((\ch{\cal P})^{\vee}p_1^*\sqrt{\td_{X_1}}p_2^*\sqrt{\td_{X_2}}
p_2^*(y)), y \in H^{ev}(X_2,{\Bbb Q}).
\end{align}
By Grothendieck Riemann-Roch theorem,
the following diagram is commutative.
\begin{equation}
\begin{CD}
{\mathbf D}(X_1) @>{{\cal F}_D}>> {\mathbf D}(X_2)\\
@V{\sqrt{\td_{X_1}}\ch}VV @VV{\sqrt{\td_{X_2}}\ch}V\\
H^{ev}(X_1,{\Bbb Q}) @>{{\cal F}_H}>> H^{ev}(X_2,{\Bbb Q})
\end{CD}
\end{equation}

\begin{lem}\label{lem:comm}
For $x \in H^{ev}(X_1,{\Bbb Z})$, $y \in H^{ev}(X_2,{\Bbb Z})$,
we get $\langle {\cal F}_H(x),y\rangle=\langle x,\widehat{{\cal F}}_H(y)\rangle$.
\end{lem}

\begin{pf}
\begin{equation}
\begin{split}
\langle {\cal F}_H(x),y\rangle &=
-\int_{X_2}(p_{2*}((\ch{\cal P})p_1^*\sqrt{\td_{X_1}}p_2^*\sqrt{\td_{X_2}}
p_1^*(x))y^{\vee}\\
&=-\int_{X_1 \times X_2}
((\ch{\cal P})p_1^*\sqrt{\td_{X_1}}p_2^*\sqrt{\td_{X_2}}p_1^*(x)p_2^*(y)^{\vee})\\
&=-\int_{X_1 \times X_2}p_1^*(x)
((\ch{\cal P})^{\vee}p_1^*\sqrt{\td_{X_1}}p_2^*\sqrt{\td_{X_2}}
p_2^*(y))^{\vee}\\
&=-\int_{X_1}x\{
p_{1*}((\ch{\cal P})^{\vee}p_1^*\sqrt{\td_{X_1}}p_2^*\sqrt{\td_{X_2}}
p_2^*(y))\}^{\vee}\\
&=\langle x,\widehat{{\cal F}}_H(y)\rangle.
\end{split}
\end{equation}
\end{pf}
\begin{lem}
For $x \in H^{ev}(X_1,{\Bbb Z})$,
${\cal F}_H(x)$ belongs to $H^{ev}(X_2,{\Bbb Z})$.
In particular ${\cal F}_H$ is an isometry of Mukai lattice.
\end{lem}
\begin{pf}
By [Y3, sect. 2], $[{\cal F}(x)]_1$ belongs to $H^2(X_2,{\Bbb Z})$.
By Lemma \ref{lem:comm}, we get
\begin{equation}
\begin{split}
\langle {\cal F}_H(x), 1 \rangle &=
\langle x,\widehat{\cal F}_H(1) \rangle \in {\Bbb Z},\\
\langle {\cal F}_H(x), \omega_2 \rangle &=
\langle x,\widehat{\cal F}_H(\omega_2) \rangle \in {\Bbb Z}.
\end{split}
\end{equation}
Hence ${\cal F}(x) \in H^{ev}(X_2,{\Bbb Z})$.
\end{pf}

Let $H$ and $H'$ be ample divisors on $X_1$ and $X_2$ respectively,
and $v \in H^{ev}(X,{\Bbb Z})$ a Mukai vector.
Let $U$ be an open subscheme of $M_H(v)$ such that
$\text{WIT}_i$ holds for $E \in U$ and $R^i p_{2*}({\cal P} \otimes p^*_1(E))$ belongs to $M_{H'}(w)$,
where $w=(-1)^i{\cal F}_H(v)$. 
We assume that $\codim_{M_H(v)}(M_H(v) \setminus U) \geq 2$
and $U \to M_{H'}(w)$ is birational.
We denote the image of $U$ by $V$.
We set $U=X_0$ and we denote projections
$U \times X_1 \times X_2 \to X_i$ and
$U \times X_1 \times X_2 \to X_i \times X_j$ by
$q_i$ and $q_{ij}$ respectively.
We also denote the projection
$U \times X_i \to U$ by $r_i$
and the projection $U \times X_i \to X_i$ by $s_i$.
Let ${\cal E}$ be a quasi-universal family of similitude 
$\rho$ on $U \times X_1$.
By the identification $U \to V$,
$R^iq_{02*}(q_{12}^*{\cal P} \otimes q_{01}^*{\cal E})$
becomes a quasi-universal family of similitude $\rho$ on $V \times X_2$.

\begin{prop}\label{prop:comm}
${\cal F}_H$ induces an isometry $v^{\perp} \to w^{\perp}$
and the following diagram is commutative.
\begin{equation}\label{eq:diagram}
\begin{CD}
v^{\perp} @>{(-1)^i{\cal F}_H}>> w^{\perp}\\
@V{\theta_v}VV @VV{\theta_w}V\\
H^2(K_H(v),{\Bbb Z}) @= H^2(K_H(w),{\Bbb Z})
\end{CD}
\end{equation}
where $K_H(v)$ is a fiber of an albanese map
$M_H(v) \to \Alb(M_H(v))$.
In particular, if $K_H(v)$ is irreducible symplectic,
then $\theta_v$ is an isometry of Hodge structures
if and only if $\theta_w$ is an isometry of Hodge structures.
\end{prop}
\begin{pf}
The first assertion follows from Lemma \ref{lem:comm}.
For $y \in w^{\perp}$, we see that
\begin{equation}
\begin{split}
\rho \theta_w(y) &=(-1)^i\left[r_{2*}(\ch((-1)^i R^iq_{02*}(q_{12*}^*({\cal P})
 \otimes q_{01}^*({\cal E})))
p_2^*(\sqrt{\td_{X_2}}y^{\vee}))\right]_1\\ 
&=(-1)^i\left[r_{2*}(q_{12}^*(\ch{\cal P})q_{01}^*(\ch{\cal E})
q_1^*( \td_{X_1})
p_2^*(\sqrt{\td_{X_2}}y^{\vee}))\right]_1\\
&=(-1)^i \left[r_{1*}q_{01*}(q_{01}^*(\ch {\cal E})q_1^*(\sqrt{ \td_{X_1}})
q_{12}^*(\ch {\cal P})p_1^*(\sqrt{ \td_{X_1}})
p_2^*(\sqrt{\td_{X_2}}y^{\vee}))\right]_1\\
&=(-1)^i\left[r_{1*}((\ch{\cal E})p_1^*(\sqrt{\td_{X_1}})
s_1^*( p_{1*}((\ch{\cal P})^{\vee}
p_1^*(\sqrt{ \td_{X_1}})p_2^*(\sqrt{\td_{X_2}}y))^{\vee}))\right]_1\\
&=(-1)^i \left[r_{1*}((\ch{\cal E})p_1^*(\sqrt{\td_{X_1}})
s_1^*(\widehat{{\cal F}}_H (y))^{\vee})\right]_1\\
&=\rho \theta_v((-1)^i \widehat{{\cal F}}_H (y)^{\vee}).
\end{split}
\end{equation}
Since $\widehat{{\cal F}}_H \circ {\cal F}_H=1_{H^{ev}(X_1,{\Bbb Z})}$,
we get \eqref{eq:diagram}.
\end{pf}

Let $\pi: X \to C$ be an elliptic K3 surface or an elliptic abelian surface.
Let $f$ be a fiber of $\pi$ and $\sigma$ is a section of $\pi$.
We set $v=r+(\sigma+k f)+a \omega \in H^{ev}(X,{\Bbb Z})$.
We shall choose a polarization $H=\sigma+n f, n \gg 0$.
By using Fourier-Mukai transformations,
Bridgeland [Br] constructed a birational map
$M_H(v) \cdots \to \Pic^0(X) \times \Hilb_X^m$,
where $2m+2=\dim M_H(v)$.
Moreover if $r \geq 3$, then this birational map
is defined by Fourier-Mukai Transformation on the complement of 
a codimension 2 subset of $M_H(v)$.  
So we can apply Proposition \ref{prop:comm}. 
By deformation arguments which are more complicated than
 those in 2.2, we can reprove Theorem \ref{thm:period}
for $r \geq 3$.

\vspace{1pc}

{\it Acknowledgement.}
I would like to thank Professor T. Katsura for valuable suggestions.
I would also like to thank Max Planck Institut f\"{u}r Mathematik
for support and hospitality.

\end{document}